\documentclass[11pt,reqno]{amsart}

\usepackage{graphicx}
\usepackage[margin=1.3in]{geometry}
\usepackage{amsmath}
\usepackage{amssymb}
\usepackage{amsthm}
\usepackage{amsfonts}
\usepackage{mathrsfs}
\usepackage[mathscr]{euscript}
\usepackage[latin1]{inputenc}
\usepackage[numbers]{natbib}
\usepackage[usenames, dvipsnames]{color}
\usepackage{enumerate}
\usepackage[scriptsize,hang,raggedright]{subfigure}
\usepackage[cmyk]{xcolor}
\usepackage{mathtools}
\usepackage{graphicx,color,amsmath,amssymb,mathrsfs}
\usepackage{amsfonts,amsthm,epsfig,epstopdf}

\usepackage{hyperref}
\usepackage{pdfsync}
\usepackage{dsfont}

\usepackage{bbm}
\usepackage{multirow}

\allowdisplaybreaks

\mathtoolsset{showonlyrefs}

\definecolor{refkey}{gray}{.45}
\definecolor{labelkey}{gray}{.45}

\definecolor{grey}{rgb}{.7,.7,.7}

\newtheorem{theorem}{Theorem}[section]
\newtheorem{proposition}[theorem]{Proposition}
\newtheorem{lemma}[theorem]{Lemma}

\theoremstyle{remark}
\newtheorem{remark}[theorem]{Remark}
\theoremstyle{definition}

\newtheorem{example}[theorem]{Example}

\def\H{\mathcal H}

\def\R{\mathbb R}

\def\e{\varepsilon}

\def\S{\Sigma}

\def\l{\lambda}

\def\pa{\partial}
\def\trace{{\rm tr}}

\def\E{\mathcal{E}}

\def\P{\mathcal{P}}

\def\Hi{\mathcal{H}}

\renewcommand{\>}{\rangle}

\newcommand{\Rn}{{\mathbb{R}^n}}

\newcommand{\ka}{{\kappa}}

\def\H{\mathcal H}

\def\l{{\lambda}}
\newcommand{\vertiii}[1]{{\left\vert\kern-0.25ex\left\vert\kern-0.25ex\left\vert #1
    \right\vert\kern-0.25ex\right\vert\kern-0.25ex\right\vert}}
\def\S{\mathbb{S}}

\def\trace{{\rm tr}}

\DeclareMathOperator{\DIVV}{div}

\newcommand{\na}{{\nabla}}
\def\>{{\rangle}}

\newcommand{\ba}{\begin{array}}
\newcommand{\ea}{\end{array}}

\newcommand{\bthm}{\begin{theorem}}
\newcommand{\ethm}{\end{theorem}}
\newcommand{\bprop}{\begin{proposition}}
\newcommand{\eprop}{\end{proposition}}
\newcommand{\blemma}{\begin{lemma}}
\newcommand{\elemma}{\end{lemma}}
\newcommand{\bexmpl}{\begin{example}}
\newcommand{\eexmpl}{\end{example}}

\newcommand{\beqn}{\begin{equation}}
\newcommand{\eeqn}{\end{equation}}
\newcommand{\beqns}{\begin{equation*}}
\newcommand{\eeqns}{\end{equation*}}

\newcommand{\pt}{\partial}

\newcommand{\Hn}{\mathcal{H}^{n-1}}

\newcommand{\V}{\mathcal{V}}

\renewcommand{\leq}{\leqslant}
\renewcommand{\geq}{\geqslant}

\definecolor{mygreen}{rgb}{0.1,0.75,0.2}

\newcounter{myenumi}
\DeclareMathOperator{\dive}{div}

\DeclareMathOperator*{\argmax}{arg\,max}

\DeclareMathOperator{\arccot}{arccot}
\DeclareMathOperator{\bary}{bar}

\numberwithin{equation}{section}

\title[On Minimizers of an anisotropic liquid drop model]{On minimizers of an anisotropic liquid drop model}

\author{Oleksandr Misiats}
\address{\parbox{\linewidth}{Department of Mathematics and Applied Mathematics,\\
Virginia Commonwealth University, Richmond, VA}}
\email{omisiats@vcu.edu}

\author{Ihsan Topaloglu}
\address{\parbox{\linewidth}{Department of Mathematics and Applied Mathematics,\\
Virginia Commonwealth University, Richmond, VA}}
\email{iatopaloglu@vcu.edu}

\date{\today}                                        

\begin{document}

\begin{abstract}
We consider a variant of Gamow's liquid drop model with an anisotropic surface energy. Under suitable regularity and ellipticity assumptions on the surface tension, Wulff shapes are minimizers in this problem if and only if the surface energy is isotropic. We show that for smooth anisotropies, in the small nonlocality regime, minimizers converge to the Wulff shape in $C^1$-norm and quantify the rate of convergence. We also obtain a quantitative expansion of the energy of any minimizer around the energy of a Wulff shape yielding a geometric stability result. For certain crystalline surface tensions we can determine the global minimizer and obtain its exact energy expansion in terms of the nonlocality parameter.
\end{abstract}

\maketitle

\baselineskip=13pt

\section{Introduction}\label{sec:intro}

In this paper we consider an anisotropic nonlocal isoperimetric problem given by
		\beqn \label{eqn:aniso_drop_rescaled}
				\inf\, \Big\{\E_\gamma(F)  \,\,  \Big\vert\,    \ |F| =1 \Big\}
		\eeqn
over sets of finite perimeter $F\subset \R^n$ where
		\beqn \label{eqn:energy}
			\E_\gamma(F) := \int_{\pa^*F}f(\nu_F)\,d\H^{n-1} + \gamma\, \int_F \! \int_F \frac{1}{|x-y|^\alpha}\,dxdy
		\eeqn 
with $0<\alpha<n$ and $| \cdot |$ denotes the Lebesgue measure.

The first term in $\E_\gamma$ is the anisotropic surface energy
	\[
		\P_f(F):=\int_{\pa^* F}f(\nu_F)\,d\H^{n-1}
	\]
defined via a one-homogeneous and convex surface tension $f:\R^n\to[0,\infty)$ that is positive on $\R^{n}\setminus\{0\}$. Here $\H^{n-1}$ is the $(n-1)$-dimensional Hausdorff measure, and $\pt^* F$ denotes the reduced boundary of $F$, which consists of points $x\in\pt F$ where the limit $\nu _{F}(x)=\lim _{\rho \to 0}\frac{-\nabla\chi _{F}(B_{\rho }(x))}{|\nabla\chi _{F}|(B_{\rho }(x))}$
exists and has length one. This limit is called the measure-theoretic outer unit normal of $F$.

The second term in the energy $\E_\gamma$ is given by the Riesz interactions
	\[
		\V(F):=\int_F\!\int_F \frac{1}{|x-y|^\alpha}\,dxdy
	\]
for $0<\alpha<n$.

The minimization problem \eqref{eqn:aniso_drop_rescaled} is equivalent (via the rescaling $\gamma=m^{(n+1-\alpha)/n}$) to the anisotropic liquid drop model
	\begin{equation}\label{eqn:aniso_drop}
		\inf\, \Big\{\mathcal{E}(E):= \P_f(E)+ \V(E) \,\,  \Big\vert\,    \ |E| =m \Big\}
	\end{equation}
introduced by Choksi, Neumayer and the second author in \cite{ChNeTo2018} as an extension of the classical liquid drop model.

Gamow's liquid drop model, initially developed to predict the mass defect curve and the shape of atomic nuclei, dates back to 1930 \cite{Ga1930}; however, it recently has generated considerable interest in the calculus of variations community (see e.g. \cite{AlBrChTo2017_3,BoCr14,ChPe2010,FrLi2015,Ju2014,KnMu2014,KnMuNo2016,LuOtto2014} as well as \cite{ChMuTo2017} for a review). The version of this model in the language of the calculus of variations includes two competing forces: 
an  attractive isotropic surface energy associated with the depletion of nucleon
density near the nucleus boundary,  and a repulsive Coulomb energy due to the interactions of
positively charged protons. These two forces are in direct competition. The surface energy prefers uniform, symmetric and connected domains whereas the repulsive term is minimized by a sequence of sets diverging infinitely apart. The parameter of the problem ($\gamma$ in \eqref{eqn:aniso_drop_rescaled} or $m$ in \eqref{eqn:aniso_drop}) sets a length scale between these competing forces.  As such, the liquid drop model is a paradigm for shape optimization via competitions of short- and long-range interactions and it appears in many different systems at all length scales.

In the anisotropic extension of the liquid drop model the global minimizer of the surface energy $\P_f(E)$ over sets $|E|=m$ is (a dilation or translation of) the {\it Wulff shape} $K_f$ associated with $f$ (cf. \cite{brothersmorgan,fonseca_wulff_rev,fonsecamuller_wulff}),
where 
\begin{equation}\label{Wulffshape}
  K_f\, : =\, \bigcap_{\nu\in\S^{n-1}}\big\{x\in\R^{n}  \, \big\vert\,  x\cdot\nu<f(\nu)\big\}.
\end{equation}
Properties of $K_f$ depend on the regularity of the surface tension $f$.

In the literature two important classes of surface tensions are considered:
	\begin{itemize}
		\item We say that $f$ is a {\it smooth elliptic surface tension} if $f \in C^\infty(\R^n \setminus\{0\})$ and there exist constants $0<\l\leq\Lambda<\infty$ such that for every $\nu\in\S^{n-1}$,
			\[
 				 \l\,|\tau|^2\leq \nabla^2f(\nu)[\tau,\tau]\leq\Lambda\,|\tau|^2\,
			\]
for all $\tau \in \R^n$ with $\tau\cdot \nu =0.$
For such surface tensions, the corresponding Wulff shape has $C^\infty$ boundary and is uniformly convex.

		\item  We say that $f$ is a {\it crystalline} surface tension if for some $N$ finite and $x_i \in \R^n,$
			\[
				f(\nu) = \max_{1\leq i\leq N} x_i \cdot \nu.
			\]
For crystalline surface tensions, the corresponding Wulff shape $K$ is a convex polyhedron.
	\end{itemize}
	
In the anisotropic liquid drop model \eqref{eqn:aniso_drop} the competition which leads to an energy-driven pattern formation is not only between the attractive and repulsive forces, as they scale differently in terms of the mass $m$, but there is also a competition between the anisotropy in the surface energy and the isotropy in the Coulomb-like energy. As shown in \cite[Theorem 3.1]{ChNeTo2018}, the problem \eqref{eqn:aniso_drop} admits a minimizer when $m$ is sufficiently small and fails to have minimizers for large values of $m$. However, \cite[Theorem 1.1]{ChNeTo2018} shows that when $f$ is smooth the Wulff shape $K_f$ is not a critical point of the energy $\E$ for any $m>0$. On the other hand, for particular crystalline surface tensions the authors prove that the corresponding Wulff shape is the unique (modulo translations) minimizer for sufficiently small $m$.
This demonstrates a fundamentally interesting situation: the regularity and ellipticity of the surface tension $f$ determines whether the isoperimetric set $K_f$ is also a minimizer of the perturbed problem \eqref{eqn:aniso_drop_rescaled}. As stated in \cite{ChNeTo2018}, while the regularity and ellipticity of the surface tension affect typically quantitative aspects of anisotropic isoperimetric problems, here, due to the incompatibility of the Wulff shape with the Riesz energies, qualitative aspects of the problem are effected, too.

Motivated by the results in \cite{ChNeTo2018}, we study qualitative properties of the minimizers of \eqref{eqn:aniso_drop_rescaled} for smooth anisotropies in the asymptotic regime $\gamma\to 0$, and obtain
	\begin{itemize}
		\item the convergence of the minimizers to the Wulff shape in strong norms, providing the rate of convergence, and
		\item an expansion of the energy around the energy of a Wulff shape in terms of $\gamma$.
	\end{itemize}
	
In particular, our first main result shows that the minimizers of $\E_\gamma$ are close to the Wulff shape in $C^1$-norm in the small $\gamma$ regime. Further we obtain quantitative estimates on how much a minimizer $F$ of $\E_\gamma$ differs from the Wulff shape when $\gamma$ is sufficiently small.

\medskip

	\bthm \label{thm:characterization_thm}
		Let $f$ be a smooth elliptic surface tension and $F$ be a minimizer of the problem \eqref{eqn:aniso_drop_rescaled}. Let $K$ denote the Wulff shape corresponding to $f$ rescaled so that $|K|=1$. Then we have the following two statements.
			\begin{itemize}
					\item[(i)] 	For $\gamma>0$ sufficiently small there exists $\psi \in C^1(\pt K)$ such that
										\[
											\pt F = \big\{ x+\psi(x)\nu_{K}(x)  \, \big\vert \, x\in \pt K \big\} 
										\]
									and
										\[
											 |F \triangle K| \lesssim \| \psi \|_{C^1(\pt K)} \lesssim |F \triangle K|^{1/(n+1)}.
										\]
						\vspace{0.1cm}
					\item[(ii)] For $\gamma>0$ sufficiently small, we have that 
										\[
											|F \triangle K| \simeq \gamma.
										\]		 
		\end{itemize}
	\ethm
	
\medskip
	
Combining parts (i) and (ii) of the theorem, we conclude that
	\[
		\gamma \lesssim \| \psi \|_{C^1(\pt K)} \lesssim \gamma^{1/(n+1)}.
	\] 	
Using the elliptic regularity theory, via Schauder estimates on the Euler--Lagrange equation, implies that $\psi \to 0$ in some $C^{2,\beta}$-norm as $\gamma\to 0$; however, finding the rate of convergence explicitly in terms of $\gamma$ seems to be a challenging task. As for quantifying the convergence rate in stronger norms, adapting arguments from Figalli and Maggi's work on the shapes of liquid drops (cf. \cite{FiMa2011}), it is possible to obtain quantitative convexity estimates on minimizers $F$ ultimately yielding $C^2$-control on the function $\psi$ via an upper bound that depends on $\gamma$ (see Remark \ref{rem:convex}). Although the result above only establishes an explicit $C^1$-control of the function $\psi$, our proof relies only on a simple geometric argument we present in the next section.

\bigskip

Next we show that the energy difference between a minimizer and the Wulff shape scales as $\gamma^2$.

	\bthm \label{thm:energy_expan}
	 Suppose $f$ is a smooth elliptic surface tension that is not a constant multiple of the Euclidean distance. Let $F$ be a minimizer of the energy $\E_\gamma$. Then for $\gamma$ sufficiently small,
	 	\beqn	\label{eqn:energy_expan}
	 		\E_\gamma(K) - \E_\gamma(F) \simeq \gamma^2
	 	\eeqn 
	where $K$ is the Wulff shape corresponding to $f$ rescaled so that $|K|=1$ and translated to have the same barycenter as $F$. 
	\ethm

\medskip

 Combined with the estimate on the symmetric difference, this expansion also yields a geometric stability estimate of the form
	\[
		\E_\gamma(K) - \E_\gamma(F) \geq C|F\triangle K|^2
	\]
for the minimizer $F$ in the small $\gamma$ regime. We prove Theorems \ref{thm:characterization_thm} and \ref{thm:energy_expan} in Section \ref{sec:characterization}.

In two dimensions, when the Wulff shape is given by a particular perturbation of a set that is symmetric with respect to the coordinate axes and lines $y=\pm x$, it is possible to determine the constant in the lower bound $\E_\gamma(K_f)-\E_\gamma(F) \geq C\,\gamma^2$, explicitly. This result is independent of the regularity of the surface tension $f$ and applies to both smooth and crystalline cases. Furthermore, when the surface tension is given by $f(\nu)=\frac{1}{2}\big(a_0|\nu\cdot e_1|+a_0^{-1}|\nu\cdot e_2|\big)$ for some $a_0>1$, we show that the minimizer of $\E_\gamma$ is a rectangle with dimensions determined explicitly in terms of $a_0$ and $\gamma$, and we obtain an expansion of the energy $\E_\gamma$ of a minimizer in terms of $\gamma$ and $a_0$ only. We prove these results in Section \ref{sec:2d}.

Finally, we would like to note that a similar incompatibility occurs also in a nonlocal isoperimetric problem considered by Cicalese and Spadaro \cite{cicalesespadaro} where the authors study the isotropic version of the energy $\E_\gamma$ (i.e., with $f$ given by the Euclidean distance) on a bounded domain $\Omega$. Here the incompatibility is due to the boundary effects. As a result of the boundary effects the isoperimetric region (in this case a ball) is not a critical point of the nonlocal term, and the authors study the asymptotic properties of the minimizers in the small $\gamma$ limit.

\subsection*{Notation} Throughout the paper  we use the notation $f\lesssim g$ to denote that $f \leq C g$ for some constant $C>0$ independent of $f$. We also write $ f\simeq g$ to denote that $c\, g \leq f \leq C\, g$ for constants $c,\,C>0$ independent of $f$. The constants $C$ we use might change from line to line unless defined explicitly. Also, when necessary, we emphasize the dependence of the constants to the parameters. In order to simplify notation we will denote the Wulff shape by $K$, suppressing the dependence on the surface tension $f$.

\section{Proofs of Theorem \ref{thm:characterization_thm} and Theorem \ref{thm:energy_expan}}\label{sec:characterization}

The proof of the first part of Theorem \ref{thm:characterization_thm} relies on a result which is independent of the optimality of the set $F$, and is rather a general property of two sets where the boundary of one of the sets is expressed as a graph over the boundary of the other set. We state this geometric result as a separate lemma since it might be of interest to readers beyond its connection to the anisotropic liquid drop model.

	\blemma \label{lem:golden_corral}
		Suppose $E$ and $F$ are bounded subsets of $\R^n$ with $C^1$ boundaries such that 
			\[
				\pt F = \big\{ x+\psi(x)\nu_{E}(x)  \, \big\vert \, x\in \pt E \big\}
			\]
		for some function $\psi \in C^1(\pt E)$.  
		\begin{itemize}
		\item[(i)] Then
		\[
	   \frac{1}{\Hn(\pt E)} |E \triangle F| \leq \| \psi \|_{C^1(\pt E)}.
		\]
		\item[(ii)] If, in addition, $\pt E \in C^2$ and  $F$ is convex, then there exits a constant $C>0$, depending only on the maximal principle curvature $\ka$ and the perimeter of $E$ (see \eqref{C} for the explicit dependence) such that 
		\begin{equation}\label{GC2}
		\| \psi \|_{C^1(\pt E)} \leq C |E \triangle F|^{\frac{1}{n+1}}.
		\end{equation}
		\end{itemize}
	\elemma

\begin{remark} In fact, the simple estimate below shows that the symmetric difference can be controlled by the $C^0$-norm of $\psi$. Namely, $|E \triangle F| \leq \Hn(\pt E)\,\|\psi\|_{C^0(\pt E)}$.
\end{remark}

	\begin{proof}
The proof of (i) is straightforward:
\[
 |E \triangle F| = \int_{\pt E}|\psi(x)| \, d\Hn \leq \|\psi\|_{C^1(\pt E)} \Hn(\pt E).
\]
To show (ii), let $z^*=(z_1^*, ..., z_n^*)= \argmax_{\pt E} |\nabla \psi(x)|$. Suppose $z^* \in \pt E \cap \pt F$.
Further, assume that in some neighborhood $U_{z^*} \subset \R^n$ there exist $f,\, g: \R^{n-1} \to \R$ such that 
\[
	\pt F = \big\{x = (x_1,...,x_n) \in U_{z^*} \, \big\vert \,  x_n = f(x_1,...,x_{n-1})\big\}
\]
and 
\[
	\pt E = \big\{x = (x_1,...,x_n) \in U_{z^*} \, \big\vert \, x_n = g(x_1,...,x_{n-1})\big\}.
\]
Since the mean curvature of $\pt E$ is bounded, we can choose the neighborhood $U_{z^*}$ such that $|U_{z^*}| \geq c_E>0$ for some constant $c_E$ depending only on $\ka$, the maximal principle curvature of $\pt E$.
Finally, without loss of generality, we may choose an appropriate rotation and translation to have $z^* = 0$,
\beqn \label{eqn:g_and_f}
g(0) = f(0) = 0 \qquad \text{ and } \qquad \nabla f(0)  = 0.
\eeqn
In this case 
\[
\max_{\pt E}| \nabla \psi(x)| = |\nabla g(0)|.
\]
Denote $\tilde{x} := (x_1,...,x_{n-1})$ and $\tilde{U}_0=\big\{\tilde{x}\in \R^{n-1}\colon (\tilde{x},0)\in U_0\big\}$. If we expand $g$ in the neighborhood $\tilde{U}_0$ of $0$, we get
\[
g(\tilde{x}) = \nabla g(0) \cdot \tilde{x} +  \frac{1}{2} \tilde{x} \cdot \nabla^2 g(0) \cdot {\tilde{x}}^{T} + o(|\tilde{x}|^2)
\]
Furthermore $\tilde{x} \cdot \nabla^2 g(0) \cdot {\tilde{x}}^{T} \geq - c_{\ka}^2 |\tilde{x}|^2$, where $c_{\ka}>0$ depends only on $\kappa$. Hence in $\tilde{U}_0$ we have
\[
g(\tilde{x}) \geq  \nabla g(0) \cdot \tilde{x} - c_{\kappa}^2 |\tilde{x}|^2 = \frac{|\nabla g(0)|^2}{4 c_{\kappa}^2} - \left|c_{\kappa} \tilde{x} - \frac{\nabla g(0)}{2 c_{\kappa}}\right|^2 
\]
On the other hand, due to the convexity of $F$, we have $f(\tilde{x}) \leq 0$ in $\tilde{U}_{0}$. Therefore,
	\begin{align*}
		|E \triangle F| &\geq \int_{\tilde{U}_{0}}|g(\tilde{x}) - f(\tilde{x})| \, d \Hn_{\tilde{x}} \geq \int_{\tilde{U}_{0} \cap \{g(\tilde{x}) \geq 0\}} g(\tilde{x})\, d\Hn_{\tilde{x}} \\
						  &\geq \int_{\tilde{U}_0 \cap \left\{\left|c_{\kappa} \tilde{x} - \frac{\nabla g(0)}{2 c_{\kappa}}\right|^2 \leq  \frac{|\nabla g(0)|^2}{4c_{\kappa}^2}\right\}}  \left(\frac{|\nabla g(0)|^2}{4 c_{\kappa}^2} - \left|c_{\kappa} \tilde{x} - \frac{\nabla g(0)}{2c_{\kappa}}\right|^2 \right) \, d\Hn_{\tilde{x}} \\
						  &=\frac{1}{c_{\kappa}^{n-1}} 
\int_{\tilde{U}_0\cap\left\{|\tilde{y}|^2 \leq  \frac{|\nabla g(0)|^2}{4c_{\kappa}^2}\right\}}  \left(\frac{|\nabla g(0)|^2}{4 c_{\kappa}^2} - \left|\tilde{y}\right|^2 \right) \,d \Hn_{\tilde{y}} \\
						  &= \frac{c_E\,\omega_{n-1}}{c_{\kappa}^{n-1}} \left(\frac{|\nabla g(0)|}{2c_{\kappa}}\right)^{n+1} - (n-1) \frac{c_E\omega_{n-1}}{c_{\kappa}^{n-1}} \int_{0}^{\frac{|\nabla g(0)|}{2c_{\kappa}}} r^{n} \,dr  \\							  
						  &= \frac{c_E\,\omega_{n-1}|\nabla g(0)|^{n+1}}{(n+1) 2^n c_{\kappa}^{2n}}.
	\end{align*}
If $z^* \in \pt E \setminus \pt F$, on the other hand, we may again choose an appropriate rotation and translation so that $z^*=0$; however, now the functions $f$ and $g$ in \eqref{eqn:g_and_f} satisfy
	\[
		g(0)=c_g, \quad f(0)=0 \qquad \text{ and } \nabla f(0)=0.
	\]
This, in turn, implies that $g(\tilde{x}) \geq c_g + \frac{|\nabla g(0)|^2}{4 c_{\kappa}^2} - \left|c_{\kappa} \tilde{x} - \frac{\nabla g(0)}{2 c_{\kappa}}\right|^2$, and estimating as above we get
	\[
		|E \triangle F| \geq |c_g| c_E + \frac{c_E\,\omega_{n-1}|\nabla g(0)|^{n+1}}{(n+1) 2^n c_{\kappa}^{2n}} \geq \frac{c_E\,\omega_{n-1}|\nabla g(0)|^{n+1}}{(n+1) 2^n c_{\kappa}^{2n}}.
	\]

Thus
\[
\max_{\pt E}|\nabla \psi(x)| \leq \left(\frac{2^n c_{\kappa}^{2n} (n+1)}{c_E\,\omega_{n-1}}\right)^{\frac{1}{n+1}} \left|E \triangle F\right|^\frac{1}{n+1}.
\]
Under the assumption that $\psi(0) = 0$, we have
\[
\max_{\pt E}| \psi(x)| \leq \max_{\pt E}| \nabla \psi(x)| \Hn(\pt E),
\]
hence altogether
\[
\|\psi\|_{C^1(\pt E)} \leq C \left|E \triangle F\right|^\frac{1}{n+1}
\]
where
\begin{equation}\label{C}
C :=\big(1+ \Hn(\pt E)\big) \left(\frac{2^n c_{\kappa}^{2n} (n+1)}{c_E\, \omega_{n-1}}\right)^{\frac{1}{n+1}}
\end{equation}
	is independent of $F$.
\end{proof}
	
\medskip
	
\begin{remark}
The constants $c_{\kappa}$ and $c_E$ differ from the actual maximal principle curvature of $\pt E$  by a factor, proportional to $(1+|\nabla g(0)|^2)^{3/2}$.
\end{remark}

\medskip

\begin{remark}
The inequality \eqref{GC2} holds true for any $C^2$ smooth set $F$, not necessarily convex. However, in this case, the constant $C$ in \eqref{C} will depend on the principle curvatures of both $F$ and $E$. In our case, we will apply the lemma in the situation, where the curvature of $E$ is known while the curvature of $F$ is not known a priori, hence the convexity of $F$ is crucial. 
\end{remark}	
	
\medskip	

For a smooth elliptic surface tension $f$ the first variation of $\P_f(E)$ with respect to a variation generated by $X\in C^{1}_c(\R^n,\R^n)$ is given by 
\[
\delta \P_f(E)[X] = \int_{\pa^*E} \DIVV^{\pa^* E}\big(\na f\circ \nu_E\big) X\cdot\nu_E \, d\H^{n-1}.
\]
Here $\dive^{\pa^* E}$ denotes the tangential divergence along $\pa^* E$. The function $H_E^f:\pa^* E\to \R$ defined by $H_E^f= \DIVV^{\pa^* E}\big(\na f\circ \nu_E\big)$ is called the \emph{anisotropic mean curvature} of the reduced boundary of $E$. The first variation of $\mathcal{V}(E)$ with respect to a variation generated by $X\in C^{1}_c(\R^n,\R^n)$, on the other hand, is given by 
\[
\delta \V(E)[X] = \int_{\pa^*E} v_E(x) X\cdot\nu_E \, d\Hi^{n-1}\,,
\]
with $v_E(x)=\int_E |x-y|^{-\alpha}\,dy$.

We say that a set $E$ is a \emph{critical point} of \eqref{eqn:aniso_drop_rescaled} if $\delta(\P_f(E) + \gamma\,\V(E))[X]=0$ for all variations with $\int_{\pa^*E} X\cdot \nu_E\, d \Hi^{n-1}=0$, i.e., variations that preserve volume to first order. Hence, a volume-constrained critical point $E$ of \eqref{eqn:aniso_drop_rescaled} satisfies the Euler-Lagrange equation
\begin{equation}\label{eqn: EL}
H_E^f(x) + v_E(x) = \lambda  \qquad \text{ for all }x \in \pa^* E,
\end{equation}
where the constant $\lambda$ is the Lagrange multiplier associated with the volume constraint $|E|=1$.

In order to obtain the rate of the $L^1$-convergence of the minimizing sets in terms $\gamma$, we will utilize the following lemma, which provides a lower bound on the energy deficit and is also a fundamental part of the proof of Theorem \ref{thm:energy_expan}.

\blemma \label{lem:energy_lower_bd}
	Suppose $f$ is a smooth elliptic surface tension that is not a constant multiple of the Euclidean distance. Let $F$ be a minimizer of the energy $\E_\gamma$. Then for $\gamma$ sufficiently small,
	 	\beqn	\label{eqn:energy_lower_bd}
	 		\E_\gamma(K) - \E_\gamma(F) \geq C \gamma^2
	 	\eeqn 
	where $K$ is the Wulff shape corresponding to $f$ rescaled so that $|K|=1$ and translated to have the same barycenter as $F$ and the constant $C$ depends only on $f$, $\alpha$, and $d$. 
\elemma

\begin{proof}
	First note that since $f$ is not a multiple of the Euclidean distance the corresponding Wulff shape $K$ is not a ball but since $K$ minimizes the perimeter functional its anisotropic mean curvature $H_K^f$ is constant on $\pt K$. On the other hand, characterization results \cite[Theorem 1.3]{ChNeTo2018} and \cite[Theorem 4.2.]{Gomez_et_al_2019} state that the only sets $E$ for which the Riesz potential $v_E$ is constant on $\pt E$ are given by balls. Hence, $K$ does not satisfy \eqref{eqn: EL}, and therefore it is not a critical point of the energy $\E_\gamma$ for any $\gamma$. This implies that there exists a function $\varphi: \pt K \to \R$ such that the first variation of $\V$ in the normal direction $\nu_K$ is negative for small perturbations by the function $\varphi$. That is, 
	\[
		\mu_2(K) := \delta \V(K)=[\varphi \nu_K] = \frac{d}{d\e} \V(K_{\varphi,\e}) \Big\vert_{\e=0} < 0
	\]
where $K_{\varphi,\e}:=\{ x + \e \varphi(x) \nu_K(x) \, \big\vert \, x\in K \}$.

Now, let
	\[
		\mu_1(K) := \int_{\pt K} D^2 f(\nabla\varphi,\nabla\varphi)  - \varphi^2 \trace(D^2 f A_K^2) \,d \Hn
	\]
where $A_K$ denotes the second fundamental form of $\pt K$. Then $\mu_1(K)=\delta^2 \P_f(K)[\varphi \nu_K]$, the second variation of $\P_f$ at $K$ (cf. \cite[Theorem 4.1]{ClMo2004}). Since $f$ is uniformly elliptic by assumption, using \cite[Lemma 4.1]{Ne2016} and arguing as in the proof of \cite[Proposition 1.9]{Ne2016} (where we also use that $\bary K = \bary F$ with $\bary$ denoting the barycenter of a set), we have that
	\[
		\mu_1(K) \geq C \int_{\pt K} |\nabla \varphi|^2\,d\Hn >0.
	\]

Using the minimality of $F$ and expanding the energies in terms of $\e$ we obtain
	\begin{align*}
		\E_\gamma(K) - \E_\gamma(F) &\geq \E_\gamma(K) - \E_\gamma(K_{\varphi,\e}) \\
													  &=     \P_f(K) - \P_f(K_{\varphi,\e}) + \gamma\,\big( \V(K) - \V(K_{\varphi,\e}) \big) \\
													  &=     -\mu_2(K) \gamma\,\e - \frac{\mu_1(K)}{2} \, \e^2 - \gamma \, o(\e).
	\end{align*}
Optimizing in $\e$ we let $\e=\big(-\mu_2(K)/\mu_1(K)\big)\gamma$. Then $\E_\gamma(K) - \E_\gamma(F) \geq C \gamma^2$ for some constant $C>0$; hence, we obtain the lower bound as claimed.	
\end{proof}

\medskip

Another important ingredient in the proof of the theorem is the regularity of quasiminimizers of the surface energy $\P_f$. We say that $F$ is a \emph{$q$-volume-constrained quasiminimizer} of $\P_f$ if
	\[
		\P_f(F) \leq \P_f(E) + q |F \triangle E| \qquad \text{ for all } E \text{ with } \quad |E|=|F|.
	\]
We are now ready to prove the theorem.

\bigskip

	\begin{proof}[Proof of Theorem \ref{thm:characterization_thm}]
	We start by noting that the nonlocal functional $\V$ is Lipschitz continuous with respect to the symmetric difference. To see this let $\alpha \in (0,n)$ and let $v_F:\R^n\to \R$ denote the Riesz potential of $F$ given by $v_F(x)= \int_{F} |x-y|^{-\alpha}\,dy$.	Hence, $\V(F)= \int_F v_F(x) \,dx$. Let $r=\omega_n^{-1/n}$, where $\omega_n$ denotes the volume of the unit ball in $\R^n$. Then
	\[
		\| v_F\|_{L^\infty(\R^n)} \leq \| v_{B_r(0)}\|_{L^\infty(\R^n)} = v_{B_r(0)}(0) = \frac{n\omega_n^{1-(n-\alpha)/n}}{n-\alpha}.
	\]
In fact, by \cite[Lemma 3]{Reichel09} and \cite[Proposition 2.1]{BoCr14}, $v_F$ is H\"{o}lder continuous with
\begin{equation}\label{eqn: Holder}
\|v_F\|_{C^{k,\beta}(\R^n)}\leq C(n, |F|, k, \beta)
\end{equation}
for $k=\lfloor n-\alpha \rfloor$ and $\beta\in(0,1)$ with $k+\beta <n-\alpha$.

A direct calculation shows that
\begin{align*}
\V(E) - \V(F) &= \int_{\R^n} \int_{\R^n} \frac{\chi_E(x)\chi_E(y) - \chi_F(x)\chi_F(y)}{|x-y|^\alpha } \, dx dy\\
& = \int_{\R^n} v_E(y)(\chi_E(y) - \chi_F(y))\, dy + \int_{\R^n} v_F(x)(\chi_E(x) - \chi_F(y))\, dx\\
&\leq \frac{2n\omega_n^{\alpha/n}}{n-\alpha}  |E\triangle F|.
\end{align*}
Hence, for any $E,\,F \subset \R^n$ with $|E| \leq |F|$ the functional $\V$ is Lipschitz continuous with respect to the symmetric difference with Lipschitz constant given by
	\beqn\label{eqn:cnalpha}
		c_{n,\alpha}:= \frac{2n\omega_n^{\alpha/n}}{n-\alpha}.
	\eeqn

Now, for any minimizer $F$ with $|F|=1$ of the energy $\E_\gamma$, we have that 
	\[
		\P_f(F) \leq \P_f(E) + \gamma \big( \V(E)-\V(F) \big) \leq \P_f(E) + \gamma c_{n,\alpha} |E \triangle F|
	\]
for any competitor $E$ with $|E|=1$. Thus $F$ is a $\gamma c_{n,\alpha}$-volume-constrained quasiminimizer of the surface energy $\P_f$. Classical arguments and regularity results for quasiminimizers in the literature (see e.g. \cite{alm66, SSA77, bomb82, DuzaarSteffen02}) imply that for $\gamma$ sufficiently small $\pt F$ is a $C^{2,\beta}$-hypersurface for all $\beta \in (0,\beta_0)$ with $\beta_0 := \min \{1,n-\alpha\}$. In fact, $\pt F$ can locally be written as a small $C^{2,\beta}$-graph over the boundary of the Wulff shape $K$ of mass 1, and $F$ is uniformly convex (see also \cite[Theorem 2.2]{ChNeTo2018} for a precise statement of this regularity result). Therefore, there exists $\psi\in C^1(\pt K)$ such that
	\[
				\pt F = \big\{ x+\psi(x)\nu_{K}(x)  \, \big\vert \, x\in \pt K \big\}. 
	\]
Since both $F$ and $K$ are uniformly convex, we can apply Lemma \ref{lem:golden_corral} to conclude that
	\[
				|F \triangle K| \lesssim \| \psi \|_{C^1(\pt K)} \lesssim |F \triangle K|^{1/(n+1)}.
	\]
This establishes part (i) of the theorem.

\medskip

In order to prove the second part, first we note that by the quantitative Wulff inequality \cite[Theorem 1.1]{FigalliMaggiPratelliINVENTIONES},
	\beqn	\label{eqn:quant_Wulff_ineq}
		\P_f(F) - \P_f(K) \geq C\,|F\triangle K|^2
	\eeqn
for some constant $C>0$ depending only on $n$ and $K$. Then, by minimality of $F$ and by Lipschitzianity of $\V$, we get
	\begin{align*}
		|F \triangle K|^2 &\leq C \big( \P_f(F) - \P_f(K) \big) \leq C \gamma\, \big(\V(K) - \V(F) \big) \\
							&\leq C \gamma \, |F \triangle K|.
	\end{align*}
Hence, $|F \triangle K| \leq C\gamma$, and we obtain the upper bound in part (ii).

In order to prove the lower bound, first suppose that $|\P_f(K)-\P_f(F)| \ll \gamma^2$ for any $\gamma>0$. Note that $\V(K)-\V(F) \geq C\gamma$ for some $C>0$ since otherwise $\gamma_k^{-1}\big(\V(K)-\V(F)\big) \to 0$ along a subsequence $\gamma_k$, which would imply that $\E_\gamma(K)-\E_\gamma(F) = o(\gamma^2)$ and this would contradict the estimate \eqref{eqn:energy_lower_bd}. Therefore, there exists a constant $C>0$ such that
	\[
		C \gamma^2 \leq \big(\P_f(K) - \P_f(F)\big) + \gamma\, \big(\V(K)-\V(F)\big).
	\]
Since $K$ minimizes the surface energy $\P_f$, again using the Lipschitzianity of $\V$, we can estimate the right-hand side by
	\[
		\gamma \, \big(\V(K)-\V(F)\big) \leq C \gamma\, |F \triangle K|.
	\]
Combining these two estimates yields $|F \triangle K| \geq C\gamma$.

If $\P_f(F)-\P_f(K) \geq C \gamma^2$, on the other hand, then using minimality of $F$ as above yields
	\[
		C\gamma^2 \leq \P_f(F)-\P_f(K) \leq C \gamma\, \big(\V(K) - \V(F) \big) \leq C\gamma\,|F \triangle K|.
	\]
Hence, we again obtain the lower bound $|F \triangle K| \geq C\gamma$, and combined with the upper bound this concludes the proof of the theorem.
\end{proof}

\medskip

\begin{remark}[Quantification of convexity]\label{rem:convex}
For $f\in C^{\infty}(\Rn\setminus\{0\})$ and $\alpha\in(0,n-1)$ it is possible to adapt the arguments in \cite[Theorem 2 and Remark 2]{FiMa2011} to include the nonlocal Riesz kernel $\V$ as the perturbation of $\P_f$, and obtain quantitative estimates in terms of $\gamma$ regarding the convexity of $F$. Namely, one can prove that
	\[
		\max_{\pt F} |\nabla^2 f(\nu_F) \nabla \nu_F - {\rm Id}_{T_x\pt F}| \leq C \gamma^{\frac{2n}{(n+2)(n+1-\alpha)}},
	\]
where $\nabla \nu_F$ denotes the second fundamental form of $F$. This, ultimately, provides a quantitative estimate on $\|\psi\|_{C^2(\pt F)}$ in terms of $\gamma$.
\end{remark}

\bigskip


We finish this section with an expansion of the energy of a minimizer of $\E_\gamma$ around the energy of the Wulff shape corresponding to smooth elliptic anisotropies that are not given by the Euclidean distance. The key idea here is that for such surface tensions the Wulff shape is not a critical point (in the sense of first variations by smooth perturbations) of the nonlocal part $\V$. Therefore, the energy expansion does not vanish at the first order, and contribution at order $\gamma$ is present.

\bigskip

\begin{proof}[Proof of Theorem \ref{thm:energy_expan}]
	We will prove this theorem in two parts. The upper bound follows by the Lipschitzianity of the nonlocal term $\V$ and the result of Theorem \ref{thm:characterization_thm}(ii). Namely, for sufficiently small $\gamma$ we have
	\[
		\E_\gamma(K)-\E_\gamma(F) \leq \gamma \, \big(\V(K) - \V(F) \big) \leq  \gamma \,C |F \triangle K| \leq C \gamma^2.
	\]

The lower bound, on the other hand, follows directly from Lemma \ref{lem:energy_lower_bd}.
\end{proof}

%
%

\section{Explicit Constructions in Two Dimensions}\label{sec:2d}

 For certain surface tensions in two dimensions a constant in the lower bound of the expansion \eqref{eqn:energy_expan} can be computed explicitly by a particular choice of small perturbations for which the first variation of the nonlocal part is negative. In order to determine these constants quantitatively, we consider one dimensional transformations of the Wulff shape. We will denote by $E_a$ the one-dimensional stretching of any set $E\subset\R^2$ with barycenter zero by a factor $a>0$, i.e.,
 	\beqn \label{eqn:trans}
 		E_a := \left\{ \left(\frac{x}{a},ay\right)\in\R^2 \, \Big\vert \, (x,y)\in E \right\}.
 	\eeqn
 
 Our first result gives an explicit lower bound of the energy expansion when the Wulff shape is such a transformation of a diagonally symmetric set. Examples of such symmetric sets include sets with smooth boundaries as well as regular polygons such as octagons, and they can be written as Wulff shapes of functions which possess dihedral symmetry. That is, if $D_4$ denotes the set of eight matrices in the dihedral group, then we will consider functions $f:\R^2 \to \R$ satisfying
 	\beqn \label{eqn:dihedral}
 		f(Ax) = f(x) \qquad \text{ for all } A \in D_4 \text{ and } x\in\R^2.
 	\eeqn 
 We also note that the proposition below does not make any assumptions on the regularity of the surface tension $f$, and applies to both the smooth and crystalline cases.
 
 	\bprop \label{prop:2d_lowerbd}
 		Let $f:\R^2 \to R$ be a surface tension (either smooth elliptic or crystalline) satisying \eqref{eqn:dihedral}. Let
 			\[
 				f_a(x_1,x_2) = f(ax_1, x_2/a) \qquad \text{ for any } a>0,
 			\] 
 		and let $K_{a_0}$ be the Wulff shape corresponding to $f_{a_0}$ for some $a_0>0$. Then for any minimizer $F$ of the energy $\E_\gamma$ defined via the surface tension $f_{a_0}$, and for $\gamma$ sufficiently small, we have
 			\[
 				\E_\gamma(K_{a_0})-\E_\gamma(F) \geq C\,\gamma^2
 			\]
 		where the constant $C$ is determined explicitly in terms of the second variation of $\P_{f_{a_0}}$ and the first variation of $\V$ around $K_{a_0}$.
 	\eprop
 	
 \begin{proof}
Let $K$ be the Wulff shape corresponding to the function $f$. Since $f$ satisfies \eqref{eqn:dihedral}, $K$ is symmetric with respect to the rotations and reflections in $D_4$. The symmetry of $K$ implies $\E_\gamma(K_{a_0}) = \E_\gamma(K_{1/{a_0}})$; hence, without loss of generality, we can take $a_0 > 1$. 

For any $f_a$ let $K_{f_a}$ be the corresponding Wulff shape. We claim that the set $K_{f_a}$ equals $K_a$ where $K_a$ is obtained from $K$ via the transformation \eqref{eqn:trans}. Since the sets $K_{f_a}$ and $K_a$ are convex it suffices to show that the boundary is mapped to the boundary. To see this, for any $\theta\in[-\pi,\pi]\setminus\{\pm \pi/2\}$, let $\phi = \arctan (a^2\tan\theta)$, and $\phi=\arccot(a^{-2}\cot\theta)$ if $\theta=\pm \pi/2$. Then we have
	\[
		\cos\phi = \frac{a^{-1}\cos\theta}{\sqrt{a^{-2}\cos^2\theta + a^{2}\sin^2\theta}} \qquad \text{ and }\qquad \sin\phi = \frac{a \sin\theta}{\sqrt{a^{-2}\cos^2\theta + a^{2}\sin^2\theta}}.
	\]
This yields,
	\begin{align*}
		\pt K_a &= \bigcap_{\theta \in [-\pi,\pi]} \Big\{ (x,y)\in\R^2 \colon a^{-1}x\cos\theta  + a y \sin\theta = f(\cos\theta,\sin\theta) \Big\} \\
				    &= \bigcap_{\theta \in [-\pi,\pi]}  \Big\{ (x,y)\in\R^2 \colon \frac{a^{-1}x\cos\theta}{\sqrt{a^{-2}\cos^2\theta + a^{2}\sin^2\theta}}  + \frac{a y \sin\theta}{\sqrt{a^{-2}\cos^2\theta + a^{2}\sin^2\theta}} \\
				    &\qquad\qquad\qquad\qquad\qquad = f\left(\frac{\cos\theta}{\sqrt{a^{-2}\cos^2\theta + a^{2}\sin^2\theta}},\frac{\sin\theta}{\sqrt{a^{-2}\cos^2\theta + a^{2}\sin^2\theta}}\right) \Big\} \\
				    &= \bigcap_{\phi \in [-\pi,\pi]} \Big\{ (x,y)\in\R^2 \colon x\cos\phi + y\sin\phi = f(a \cos\phi,a^{-1}\sin\phi) \Big\} \\
				    &=\bigcap_{\phi \in [-\pi,\pi]} \Big\{ (x,y)\in\R^2 \colon x\cos\phi + y\sin\phi = f_a(\cos\phi,\sin\phi) \Big\} = \pt K_{f_a}.
	\end{align*}

For any $a_0$, let $K_{a_0}$ be the Wulff shape determined by the surface tension $f_{a_0}$. Then for any $a$ close to $a_0$ the perimeter can be expanded as
\beqn \label{eqn:per_expansion}
\P_f(K_{a}) = \P_f(K_{a_0}) + \frac{\mu_1(K_{a_0})}{2}(a-a_0)^2 + O\big( (a-a_0)^3 \big)
\eeqn
where 
	\[
		\mu_1(K_{a_0}):= d^2/da^2 \P_f(K_{a}) \big\vert_{a=a_0}.
	\]
Note that $d/da \P_f(K_a) \big\vert_{a=a_0}=0$  since $K_{a_0}$ is the corresponding Wulff shape; hence, it is a critical point. Moreover, as both $K_{a_0}$ and $K_{a}$ are convex and perturbations of $K$, in two dimensions they intersect at at most four points. Hence, using at most four functions, it is possible to express $\pt K_{a}$ locally as a graph over $\pt K_{a_0}$. Therefore, arguing as in the proof of Theorem \ref{thm:energy_expan} we get that $\mu_1(K_{a_0}) > 0$.

On the other hand, expanding the nonlocal term, we get
\beqn \label{eqn:nonlocal_expansion}
\V(K_{a}) = \V(K_{a_0}) + \mu_2(K_{a_0}) (a-a_0) + O\big((a-a_0)^2\big)
\eeqn
where
	\[
		\mu_2(K_{a_0}) := \frac{d}{da} \V(K_{a})\big\vert_{a=a_0}.
	\]
In order to explicitly evaluate the first variation of the nonlocal energy with respect to these special perturbations, we introduce the change of variables $\tilde{x}_i = x_i/a$ and $\tilde{y}_i = y_i/a$ for $i=1,2$. This yields
\[
\V(K_{a}) = \int_{K} \! \int_{K} \Big(a^2(\tilde{x}_1 - \tilde{x}_2)^2 + a^{-2}(\tilde{y}_1 - \tilde{y}_2)^2\Big)^{-\frac{\alpha}{2}} \, d \tilde{x}_1 d \tilde{x}_2 d \tilde{y}_1 d \tilde{y}_2.
\]
Hence,
\[
\mu_2(K_{a_0})= -\frac{\alpha}{a_0}\int_{K}\! \int_{K} \frac{a_{0}^2(\tilde{x}_1 - \tilde{x}_2)^2 - a_0^{-2}(\tilde{y}_1 - \tilde{y}_2)^2}{ \big(a_0^2(\tilde{x}_1 - \tilde{x}_2)^2 + a_0^{-2}(\tilde{y}_1 - \tilde{y}_2)^2\big)^{1+\frac{\alpha}{2}} }\, d \tilde{x}_1 d \tilde{x}_2 d \tilde{y}_1 d \tilde{y}_2.
\]

Changing the variables once again, we get
\[
\mu_2(K_{a_0})  = -\frac{\alpha}{a_0} \int_{K_{a_0}}\! \int_{K_{a_0}} \frac{(x_1-x_2)^2 - (y_1-y_2)^2}{\big((x_1-x_2)^2 + (y_1-y_2)^2\big)^{1+\frac{\alpha}{2}}}\, dx_1 dx_2 dy_1 dy_2.
\]
This show that $\mu_2(K_{a_0})<0$ due to the fact that the deformation of $K$ into $K_{a_0}$ stretches the domain in the $x$-direction, hence increasing the first term in the integral, and at the same time shrinks it in the $y$-direction, thus decreasing the second term in the integral. 

Referring back to the expansions \eqref{eqn:per_expansion} and \eqref{eqn:nonlocal_expansion}, there exists two positive constants $C_1$ and $C_2$ such that
	\[
		\E_\gamma(K_a) \leq \E_\gamma(K_{a_0}) + \frac{1}{2}\mu_1(K_{a_0})\,(a-a_0)^2 + \gamma \mu_2(K_{a_0})\, (a-a_0) + C_1 (a-a_0)^3 + \gamma C_2 (a-a_0)^2.
	\]
Optimizing in $(a-a_0)$ we let
	\[
		a-a_0 = -\frac{\mu_2(K_{a_0})}{\mu_1(K_{a_0})}\,\gamma
	\]
and note that the coefficient is positive since $\mu_2(K_{a_0})<0$.
Then using the fact that $F$ is a minimizer, we get
	\begin{align*}
		\E_\gamma (F) & \leq \E_\gamma(K_a)  \leq  \E_\gamma(K_{a_0}) - \frac{\mu_2^2(K_{a_0})}{4\mu_1(K_{a_0})}\gamma^2 + \left( C_2 \frac{\mu_2^2(K_{a_0})}{4\mu_1^2(K_{a_0})}-C_1 \frac{\mu_2^3(K_{a_0})}{8\mu_1^3(K_{a_0})} \right)\gamma^3 \\
							  & \leq \E_\gamma(K_{a_0}) - \frac{\mu_2^2(K_{a_0})}{8\mu_1(K_{a_0})}\,\gamma^2 
	\end{align*}
for $\gamma>0$ sufficiently small. Hence,
	\[
		\E_\gamma(K_{a_0}) - \E_\gamma(F) \geq \frac{\mu_2^2(K_{a_0})}{8\mu_1(K_{a_0})} \,\gamma^2
	\]
with the constant depending only on the set $K$ (that is, on the surface tension $f$) and $a_0$.
\end{proof} 	

\medskip

\begin{remark}[More general sets] The proposition above can be stated for more general Wulff shapes which are not necessarily perturbations via \eqref{eqn:trans} of a set symmetric with respect to the dihedral group $D_4$. In fact, a sufficient condition on a set $S$ for the above proof to work is that
	\begin{equation}\label{ast}
		 \int_{S}\! \int_{S} \frac{(x_1-x_2)^2 - (y_1-y_2)^2}{\big((x_1-x_2)^2 + (y_1-y_2)^2\big)^{1+\frac{\alpha}{2}}}\, dx_1 dx_2 dy_1 dy_2 \neq 0.
	\end{equation}
Sets considered in the above proposition are $S = K_{a_0}$ where $K$ can be a disk, square, regular octagon, etc. As mentioned before, the map $K \to K_{a_0}$ for $a_0 \neq 1$ is stretching (shrinking) the set $K$ in $x$ direction while shrinking (stretching) the set in $y$ direction, and is one of the examples of the perturbation for which 
\[
\frac{d}{d_{a_0}} \V(K_{a_0})\Big\vert_{a_0=1} \neq 0.
\] 
We conjecture that one can perform a similar shrinking/stretching deformation  $K \to \tilde{K}_{a_0}$ along some direction $\nu$ such that 
\[
\frac{d}{d_{a_0}} \V(\tilde{K}_{a_0})\Big\vert_{a_0=1} \neq 0.
\]
for any $K$ different from a ball.
\end{remark}

\medskip

\begin{remark}[The constants $\mu_1(K_{a_0})$ and $\mu_2({K_{a_0}})$]
While approximate values of $\mu_1(K_{a_0})$ and $\mu_2(K_{a_0})$ can be found numerically, finding their exact values analytically is a challenging task. 
Although the perturbation of the Wulff shape is given by a simple transformation, determining the exact value of $\mu_1$ would require an explicit formula for the surface tension $f$ corresponding to $K$ in order to write $f(\nu_{K_a})$ in terms of $f(\nu_{K})$.

For the constant $\mu_2$, on the other hand, we can derive estimates in different $a_0$ regimes, using the properties of the set $K$. We list these estimates here.
	\begin{enumerate}[1.] \itemsep=1.5em
			\item For $a_0 \gg 1$, we have
\[
\mu_2(K_{a_0}) = -\frac{\alpha}{a_0^{1+\alpha}} \int_{K}\int_{K} \frac{d x_1 dx_2 dy_1 dy_2}{(x_1-x_2)^{\alpha}} + o(a_0^{-(1+\alpha)}).
\]
Hence, $\lim_{a_0\to\infty} \mu_2(K_{a_0}) = 0$.

			\item Note that
				\[
					\mu_2(K) = -\alpha \int_K \! \int_K \frac{(x_1-x_2)^2 - (y_1-y_2)^2}{\big((x_1-x_2)^2+(y_1-y_2)^2\big)^{1+\frac{\alpha}{2}}}\,dx_1 dx_2 dy_1 dy_2.
				\]
			Since both the denominator of the above integral and the set $K$ is symmetric with respect to swapping the variables $x_i$ and $y_i$, we get that 
				$\mu_2(K)=0$. Hence, $K$ is a critical point of $\V$ with respect to this special class of perturbations.
			
			\item For $a_0$ close to 1, 
\[
\mu_2(K_{a_0})  = \frac{d}{da} \mu_2(K_{a})\Big\vert_{a = 1} (a_0-1)  + O(a_0-1)^2
\]
where 
\begin{multline} \nonumber
\frac{d}{da} \mu_2(K_a)\big\vert_{a = 1}  \\
		= -2 \alpha \int_{K} \! \int_{K} \frac{-\alpha\big((x_1-x_2)^4 + (y_1-y_2)^4\big) + (4+\alpha)(x_1-x_2)^2 (y_1-y_2)^2}{\big((x_1-x_2)^2 + (y_1-y_2)^2\big)^{2+\alpha/2}} \, dx_1 dx_2 dy_1 dy_2
\end{multline}
			\item We may estimate $\mu_2(K_{a_0})$ and $\frac{d}{da} \mu(K_a)\big\vert_{a = 1}$ independently of $K$. Suppose $K$ is an arbitrary convex set which is symmetric with respect to the lines $y=\pm x$, such that $\partial K$ passes through $(2p, 0)$ for some $p>0$. Then
\[
S_{\min} \subset K \subset S_{\max}
\]
where $S_{\max} = [-2p,2p] \times [-2p,2p]$, and $S_{\min} = [-p,p] \times [-p,p]$ (see Figure \ref{fig:Wulff}).

\begin{figure}[h!]
	\begin{center}
		\includegraphics[width=0.4\linewidth]{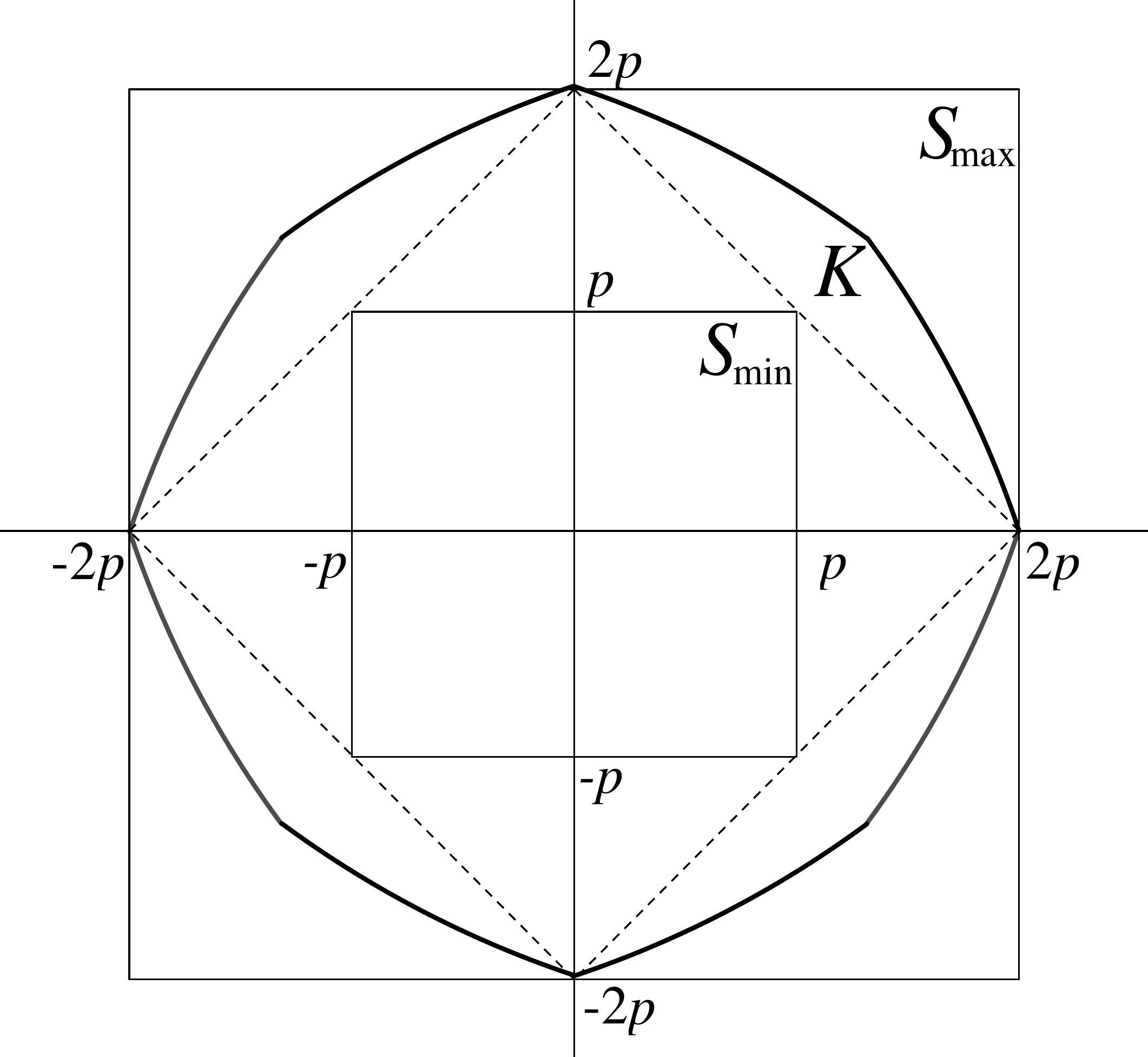}
		\caption{\footnotesize{For any Wulff shape $K$ which is convex and has 8-fold symmetry, we can find squares $S_{\min}$ and $S_{\max}$ as depicted above.}} \label{fig:Wulff}
	\end{center}
\end{figure}

Moreover,
\begin{align*}
 \int_{S_{\max}}\! \int_{S_{\max}} &\frac{a_{0}^2(x_1 - x_2)^2}{ \big(a_0^2(x_1 - x_2)^2 + a_0^{-2}(y_1 - y_2)^2\big)^{1+\frac{\alpha}{2}} }\, dx_1 dx_2 dy_1 dy_2 \\
 													&\geq \int_{K}\! \int_{K} \frac{a_{0}^2(x_1 - x_2)^2}{ \big(a_0^2(x_1 - x_2)^2 + a_0^{-2}(y_1 - y_2)^2\big)^{1+\frac{\alpha}{2}} } \,dx_1 dx_2 dy_1 d y_2  \\
													&\geq \int_{S_{\min}}\! \int_{S_{\min}} \frac{a_{0}^2(x_1 - x_2)^2}{ \big(a_0^2(x_1 - x_2)^2 + a_0^{-2}(y_1 - y_2)^2\big)^{1+\frac{\alpha}{2}} } \,dx_1 dx_2 dy_1 d y_2 \\
													&= 2^{\alpha-4} \int_{S_{\max}}\! \int_{S_{\max}} \frac{a_{0}^2(x_1 - x_2)^2}{ \big(a_0^2(x_1 - x_2)^2 + a_0^{-2}(y_1 - y_2)^2\big)^{1+\frac{\alpha}{2}} } \,dx_1 dx_2 dy_1 d y_2.
\end{align*}
Analogously,
\begin{align*}
 \int_{S_{\max}} \! \int_{S_{\max}} &\frac{a_{0}^{-2}(y_1 - y_2)^2}{ \big(a_0^2(x_1 - x_2)^2 + a_0^{-2}(y_1 - y_2)^2\big)^{1+\frac{\alpha}{2}} } \,dx_1 dx_2 dy_1 d y_2 \\
 												 & \geq \int_{K}\! \int_{K} \frac{a_{0}^{-2}(y_1 - y_2)^2}{ \big(a_0^2(x_1 - x_2)^2 + a_0^{-2}(y_1 - y_2)^2\big)^{1+\frac{\alpha}{2}} } \,dx_1 dx_2 dy_1 d y_2 \\
												 & \geq 2^{\alpha-4} \int_{S_{\max}}\! \int_{S_{\max}} \frac{a_{0}^{-2}(y_1 - y_2)^2}{ \big(a_0^2(x_1 - x_2)^2 + a_0^{-2}(y_1 - y_2)^2\big)^{1+\frac{\alpha}{2}} } \,dx_1 dx_2 dy_1 d y_2.
\end{align*}
Therefore, 
\begin{align*}
	-\frac{\alpha}{a_0} \int_{S_{\max}}\! \int_{S_{\max}} &\frac{2^{\alpha - 4} a_{0}^2(x_1 - x_2)^2 -a_{0}^{-2} (y_1 - y_2)^2}{\big(a_0^2(x_1 - x_2)^2 + a_0^{-2}(y_1 - y_2)^2\big)^{1+\frac{\alpha}{2}} } \,dx_1 dx_2 dy_1 dy_2 \geq \mu_2(K_{a_0})\\
																				&  \geq -\frac{\alpha}{a_0} \int_{S_{\max}} \! \int_{S_{\max}} \frac{a_{0}^2(x_1 - x_2)^2 -2^{\alpha - 4} a_{0}^{-2} (y_1 - y_2)^2}{\big(a_0^2(x_1 - x_2)^2 + a_0^{-2}(y_1 - y_2)^2\big)^{1+\frac{\alpha}{2}} }\, dx_1 dx_2 dy_1 dy_2 
\end{align*}
Similar upper and lower bounds can be found for $\frac{d}{da} \mu_2(K_a)\big\vert_{a = 1} $ as well. 
	\end{enumerate}
\end{remark}
 	
 \bigskip	
 
 When the Wulff shape is given by a rectangle in two dimensions, due to a rigidity theorem by Figalli and Maggi, we obtain a quantitative description of the minimizers as well as an asymptotic expansion of its energy in terms of $\gamma$ and the Wulff shape. 
 
 	\bprop \label{prop:rectangle}
 		Let $S=[-1/2,1/2]\times[-1/2,1/2]$ be the square of area 1. For $a_0>1$, let $f(\nu)=\frac{1}{2}\big(a_0|\nu\cdot e_1|+a_0^{-1}|\nu\cdot e_2|\big)$ be the surface tension whose corresponding Wulff shape is $S_{a_0}$ obtained via the transformation \eqref{eqn:trans}. Then there exists $\gamma_*>0$ such that for $\gamma<\gamma_*$ any minimizer of $\E_\gamma$ is a rectangle $S_{a}$ where
 		\beqn \label{eqn:rectangle}
 			a = a_0 - \frac{\mu_2(a_0)a_0^2}{2}\,\gamma
 		\eeqn
 	and
 		\begin{multline} \label{eqn:rectangle_energy}
 			\E_\gamma(S_a) = \E_\gamma(S_{a_0}) - \Big(\frac{\mu_2(a_0) a_0}{2}\Big)^2\, \gamma^2 \\
 					+ \Bigg( \Big(\frac{\mu_2(a_0) a_0}{2}\Big)^3 + \frac{\mu_2^2(a_0) \mu_3(a_0) a_0^4}{8} \Bigg)\, \gamma^3 + o(\gamma^3)
 		\end{multline}
 with $\mu_2(a_0)=\frac{d}{da} \V(S_{a})\big\vert_{a=a_0}$ and $\mu_3(a_0)=\frac{d^2}{da^2} \V(S_{a})\big\vert_{a=a_0}$. 
 	\eprop
  
 \begin{proof}
Let $F$ be a minimizer of $\E_\gamma$. As shown in the proof of Theorem \ref{thm:characterization_thm} above, for $\gamma$ sufficiently small, $F$ is a $\gamma c_{n,\alpha}$-volume-constrained quasiminimizer of the surface energy $\P_f$. Then, by the two dimensional rigidity theorem \cite[Theorem 7]{FiMa2011} of Figalli and Maggi, which states that if $f$ is a crystalline surface tension then any $q$-volume-constrained quasiminimizer with sufficiently small $q$ is a convex polygon with sides aligned with those of the Wulff shape, we get that $F$ is a rectangle with side parallel to $S_{a_0}$. Thus, there exists $\gamma_*>0$, such that for $\gamma<\gamma_*$ we have $F=S_{a}$ for some $a>1$.

In order to find the optimal scaling $a$, we expand the perimeter and the nonlocal term around $a_0$ and get
	\[
		\E_\gamma(S_a) = \E_\gamma(S_{a_0}) + \frac{1}{a_0^2}(a-a_0)^2 + \gamma \mu_2(a_0) (a-a_0) +\frac{\gamma}{2} \mu_3(a_0)(a-a_0)^2  -\frac{1}{a_0^3}(a-a_0)^2 + \cdots.
	\]
Optimizing at the second-order (i.e., the second and third terms in the expansion above) yields, as before,
	\[
		a-a_0 = -\frac{\mu_2(a_0) a_0^2}{2} \,\gamma,
	\]
and we obtain \eqref{eqn:rectangle}. Plugging this back into $\E_\gamma(S_a)$ we get \eqref{eqn:rectangle_energy}, i.e., an exact expansion of the energy of a minimizer in $\gamma$.
 \end{proof}
 
 \medskip
 
\begin{remark}
While we cannot determine the constant $\gamma_*$ explicitly, the expansion \eqref{eqn:rectangle_energy} yields an explicit upper bound on $\gamma_*$. Namely $\E_\gamma(S_a) < \E_\gamma(S_{a_0})$ implies that
	\[
		\gamma < \frac{2(\mu_2(a_0)a_0)^2}{(\mu_2(a_0) a_0)^3 + \mu_2^2(a_0) \mu_3(a_0) a_0^4}.
	\]
	
\end{remark}

\bigskip

\subsection*{Acknowledgments} The authors would like to thank Gian Paolo Leonardi for bringing the question of energy expansion around the energy of the Wulff shape to our attention and to Marco Bonacini, Riccardo Cristoferi and Robin Neumayer for their valuable suggestions and comments. Finally, the authors are grateful to the referees for their careful reading of the manuscript and for their detailed suggestions.


\bibliographystyle{IEEEtranS}
\def\url#1{}
\bibliography{references}

\end{document}